\documentclass{amsart}
\usepackage{amssymb}

\usepackage{graphicx}
\usepackage{amscd}

\newtheorem{theorem}{Theorem}
\theoremstyle{plain}

\newtheorem{corollary}{Corollary}

\newtheorem{definition}{Definition}
\newtheorem{example}{Example}

\newtheorem{problem}{Problem}

\newtheorem{remark}{Remark}
\input{tcilatex}

\begin{document}
\title[Subspace Homogeneous and Quotient Homogeneous Banach Spaces]{A solution of one Problem of J. Lindenstrauss and H.P. Rosenthal on Subspace
Homogeneous and Quotient Homogeneous Banach Spaces with Application to the
Approximation Problem and to the Schroeder - Bernstein Problem}
\author{Eugene V. Tokarev}
\address{B.E. Ukrecolan, 33-81, Iskrinskaya str., 61005, Kharkiv-5, UKRAINE}
\email{tokarev@univer.kharkov.ua}
\dedicatory{Dedicated to the memory of S. Banach.}

\begin{abstract}
In article is constructed a wide couple of pairwice non isomorphic separable
superreflexive Banach spaces $E$ that are subspace homogeneous. Their
conjugates are quotient homogeneous. None of this couple neither isomorphic
to its Carthesian square nor has the approximation property. At the same
time, any such $E$ is isomorphic to $E\oplus E\oplus W$ for a some Banach
space $W$ and, hence, solves the Schroeder - Bernstein problem.
\end{abstract}

\maketitle

\section{Introduction}

Notions of subspace- and of quotient-homogeneous Banach spaces were
introduced by J. Lindenstrauss and H.P. Rosenthal [1] (definitions will be
given later). They showed that $l_{2}$ and $c_{0}$ are subspace-homogeneous,
that $l_{2}$ and $l_{1}$ are quotient-homogeneous and interested in
following questions:

\begin{enumerate}
\item  \textit{Whether there exists separable subspace-homogeneous Banach
spaces that are different from }$l_{2}$\textit{\ and }$c_{0}$\textit{?}

\item  \textit{Whether there exists separable quotient homogeneous Banach
spaces that are different from }$l_{2}$\textit{\ and }$l_{1}$\textit{?}

\item  \textit{Whether there exists non separable subspace- (or quotient-)
homogeneous Banach spaces}?
\end{enumerate}

In the article it will be shown that answers on the first and on the second
question are affirmative. Any separable superreflexive Banach space $X$ may
be isomorphically (even $(1+\varepsilon )$-isomorphicaly) embedded in a
subspace-homogeneous separable superreflexive Banach space $E_{X}$ of the
same type and cotype as $X$. Of course, $\left( E_{X}\right) ^{\ast }$ is
quotient homogeneous (by superreflexivity of $E_{X}$ and by duality).
Methods that are presented below give no a chance to establish the existence
of nonreflexive B-convex spaces with such properties.

The third question is more complicated. A succeed in resolving this problem
depends on set-theoretical hypothesis. It may be proved that if there exists
a regular cardinal $\varkappa $ with a property: $\varkappa =\varkappa
^{\ast }$, where $\varkappa ^{\ast }=\sup \{2^{\tau }:\tau <\varkappa \}$,
then there exists a Banach space of dimension $\varkappa $ which is
subspace-homogeneous. Recall that the existence of a such cardinal $%
\varkappa $ cannot be proved in ZFC. Author expect to present this result in
other paper.

Unexpectedly it turned out, that spaces $E_{X}$ mentioned above posses
additional peculiar property. Namely, if there exists a real number $p\neq 2$
such that $l_{p}$ is finite representable in $X$, then $E_{X}$ does not
enjoy the approximation property.

This gives an independent solution of the \textit{approximation problem},
that was solved be P. Enflo [2].

If $X$ is not isomorphic to a Hilbert space, then $E_{X}$ and $E_{X}\oplus
E_{X}$ are not isomorphic (notice, this fact gives an independent solution
of \textit{Banach's problem on squares}, that was solved firstly in [3] and
[4]).

Any $E_{X}$ is universal for all separable spaces that are finitely
representable in $E_{X}$. This result is related to a \textit{problem on
existing of separable envelopes}: any $E_{X}$ is an ''almost envelope'': if $%
Z$ is finitely representable in $E_{X}$ then, for any $\varepsilon >0$, $Z$
is $(1+\varepsilon )$-isomorphic to a subspace of $E_{X}$.

Any $E_{X}$ has a complement in any Banach space $Y$ that contains $E_{X}$
and is finitely representable in $E_{X}$. In particular, $E_{X}\oplus E_{X}$
is isomorphic to a complemented subspace of $E_{X}$. Since $E_{X}$ obviously
is a complemented subspace of $E_{X}\oplus E_{X}$, a pair $\left\langle
E_{X},E_{X}\oplus E_{X}\right\rangle $ solves the \textit{%
Schroeder-Bernstein problem} that in other way was solved by W.T. Gowers [5].

\section{Definitions and notations}

\begin{definition}
A Banach space $X$ is said to be subspace homogeneous if for any pair of its
subspaces $Y$ and $Z$ both of equal codimension (finite or infinite) such
that there exists an isomorphism $u:Y\rightarrow Z$ there exists an
isomorphical automorphism $U:X\rightarrow X$ such that its restriction $%
U\mid _{Y}=u$.
\end{definition}

\begin{definition}
A Banach space $X$ is said to be quotient homogeneous if for any pair of its
quotients $X/Y$ and $X/Z$ such that both $Y$ and $Z$ are of equal (finite or
infinite) dimension and there exists an isomorphism $u:X/Y\rightarrow X/Z$,
there exists also an automorphism $U:X\rightarrow X$ such that $u\circ
h_{Y}=h_{Z}\circ U$,.where $h_{Y}:X\rightarrow X/Y$ and $h_{z}:X\rightarrow
X/Z$ are standard quotient maps.
\end{definition}

\begin{definition}
Let $X$, $Y$ are Banach spaces. $X$ is \textit{finitely representable} in $Y$
(in symbols: $X<_{f}Y$) if for each $\varepsilon >0$ and for every finite
dimensional subspace $A$ of $X$ there exists a subspace $B$ of $Y$ and an
isomorphism $u:A\rightarrow B$ such that $\left\| u\right\| \left\|
u^{-1}\right\| \leq 1+\varepsilon $.

Spaces $X$ and$\ Y$ are said to be finitely equivalent if $X<_{f}Y$ and $%
Y<_{f}X$.

Any Banach space $X$ generates a class 
\begin{equation*}
X^{f}=\{Y\in \mathcal{B}:X\sim _{f}Y\}
\end{equation*}
\end{definition}

Let $\mathcal{K}$ be a class of Banach spaces; $X$ be a Banach space. A
notion $X<_{f}\mathcal{K}$ (or, equivalently, $X\in \mathcal{K}^{<f}$ means
that $X<_{f}Y$ for every $Y\in \mathcal{K}$.

For any two Banach spaces $X$, $Y$ their \textit{Banach-Mazur distance }is
given by 
\begin{equation*}
d(X,Y)=\inf \{\left\| u\right\| \left\| u^{-1}\right\| :u:X\rightarrow Y\},
\end{equation*}
where $u$ runs all isomorphisms between $X$ and $Y$ and is assumed, as
usual, that $\inf \varnothing =\infty $.

It is well known that $\log d(X,Y)$ defines a metric on each class of
isomorphic Banach spaces. A set $\frak{M}_{n}$ of all $n$-dimensional Banach
spaces, equipped with this metric, is a compact metric space that is called 
\textit{the Minkowski compact} $\frak{M}_{n}$.

A disjoint union $\cup \{\frak{M}_{n}:n<\infty \}=\frak{M}$ is a separable
metric space, which is called the \textit{Minkowski space}.

Consider a Banach space $X$. Let $H\left( X\right) $ be a set of all its 
\textit{different} finite dimensional subspaces (\textit{isometric finite
dimensional subspaces of }$X$\textit{\ in }$H\left( X\right) $\textit{\ are
identified}). Thus, $H\left( X\right) $ may be regarded as a subset of $%
\frak{M}$, equipped with the restriction of the metric topology of $\frak{M}$%
.

Of course, $H\left( X\right) $ need not to be a closed subset of $\frak{M}$.
Its closure in $\frak{M}$ will be denoted $\overline{H\left( X\right) }$.
From definitions follows that $X<_{f}Y$ if and only if $\overline{H\left(
X\right) }\subseteq \overline{H\left( Y\right) }$. Certainly, spaces $X$ and 
$Y$ are \textit{finitely equivalent }( $X\sim _{f}Y$) if and only if $%
\overline{H\left( X\right) }=\overline{H\left( Y\right) }$.

Thus, there is a one to one correspondence between classes of finite
equivalence $X^{f}=\{Y\in \mathcal{B}:X\sim _{f}Y\}$ and closed subsets of $%
\frak{M}$ of kind $\overline{H\left( X\right) }$.

Indeed, all spaces $Y$ from $X^{f}$ have the same set $\overline{H\left(
X\right) }$. This set, uniquely determined by $X$ (or, equivalently, by $%
X^{f}$), will be denoted by $\frak{M}(X^{f})$ and will be referred as to 
\textit{the Minkowski's base of the class} $X^{f}$.

Using this correspondence, it may be defined a set $f\left( \mathcal{B}%
\right) $ of all different classes of finite equivalence, assuming (to
exclude contradictions with the set theory) that members of $f\left( 
\mathcal{B}\right) $ are sets $\frak{M}(X^{f})$. For simplicity it may be
says that members of $f\left( \mathcal{B}\right) $ are classes $X^{f}$
itself.

Clearly, $f\left( \mathcal{B}\right) $ is partially ordered by the relation $%
\frak{M}(X^{f})\subseteq \frak{M(}Y^{f})$, which may be replaced by the
relation $X^{f}<_{f}Y^{f}$ of the same meaning. The minimal (with respect to
this order) element of $\ f\left( \mathcal{B}\right) $ is the class $\left(
l_{2}\right) ^{f}$ (the Dvoretzki theorem); the maximal one - the class $%
\left( l_{\infty }\right) ^{f}$ (an easy consequence of the Hahn-Banach
theorem). Other $l_{p}$'s are used in the classifications of Banach spaces,
which was proposed by L. Schwartz [6].

For a Banach space $X$ its $l_{p}$-\textit{spectrum }$S(X)$ is given by 
\begin{equation*}
S(X)=\{p\in\lbrack0,\infty]:l_{p}<_{f}X\}.
\end{equation*}

Certainly, if $X\sim_{f}Y$ then $S(X)=S(Y)$. Thus, the $l_{p}$-spectrum $%
S(X) $ may be regarded as a property of the whole class $X^{f}$. So,
notations like $S(X^{f})$ are of obvious meaning.

Let $X$ be a Banach space. It is called:

\begin{itemize}
\item  $c$-\textit{convex,} if $\infty \notin S(X)$;

\item  $B$-\textit{convex,} if $1\notin S\left( X\right) $;

\item  \textit{Finite universal,} if $\infty \in S(X)$.

\item  \textit{Superreflexive,} if every space of the class $X^{f}$ is
reflexive.
\end{itemize}

Equivalently, $X$ is superreflexive if any $Y<_{f}X$ is reflexive. Clearly,
any superreflexive Banach space is $B$-convex.

\section{Quotient closed and divisible classes of finite equivalence}

In this section it will be shown how to enlarge a Minkowski's base $\frak{M}%
(X^{f})$ of a certain $B$-convex class $X^{f}$ to obtain a set $\frak{N}$,
which will be a Minkowski's base $\frak{M}(W^{f})$ for some class $W^{f}$
that holds the $B$-convexity (and, also, the $l_{p}$-spectrum) of the
corresponding class $X^{f}$ and will be quotient closed.

\begin{definition}
A class $X^{f}$ (and its Minkowski's base $\frak{M}(X^{f})$) is said to be
divisible if some (equivalently, any) space $Z\in X^{f}$ is finitely
representable in any its subspace of finite codimension.
\end{definition}

\begin{definition}
Let $\{X_{i}:i\in I\}$ be a collection of Banach spaces. A space 
\begin{equation*}
l_{2}\left( X_{i},I\right) =\left( \sum \oplus \{X_{i}:i\in I\}\right) _{2}
\end{equation*}
is a Banach space of all families $\{x_{i}\in X_{i}:i\in I\}=\frak{x}$, with
a finite norm 
\begin{equation*}
\left\| \frak{x}\right\| _{2}=\sup \{(\sum \{\left\| x_{i}\right\|
_{X_{i}}^{2}:i\in I_{0}\})^{1/2}:I_{0}\subset I;\text{ }card\left(
I_{0}\right) <\infty \}.
\end{equation*}
\end{definition}

\begin{example}
Any Banach space $X$ may be isometricaly embedded into a space 
\begin{equation*}
l_{2}(X)=(\sum\nolimits_{i<\infty }\oplus X_{i})_{2},
\end{equation*}
where all $X_{i}$'s are isometric to $X$. Immediately, $l_{2}(X)$ generates
a divisible class $\mathsf{D}_{2}(X^{f})=\left( l_{2}(X)\right) ^{f}$ which
with the same $l_{p}$-spectrum as $X^{f}$ and is superreflexive if and only
if $X^{f}$ is superreflexive too.
\end{example}

\begin{remark}
The procedure $\mathsf{D}_{2}:X^{f}\rightarrow \left( l_{2}(X)\right) ^{f}$
may be regarded as a closure operator on the partially ordered set $f\left( 
\mathcal{B}\right) $. Indeed, it is

\begin{itemize}
\item  Monotone, i.e., $X^{f}<_{f}\mathsf{D}_{2}(X^{f})$;

\item  Idempotent, i.e., $\mathsf{D}_{2}(X^{f})=\mathsf{D}_{2}(\mathsf{D}%
_{2}(X^{f}))$;

\item  Preserve the order: $X^{f}<_{f}Y^{f}$ $\Longrightarrow \mathsf{D}%
_{2}(X^{f})<_{f}\mathsf{D}_{2}(Y^{f})$.
\end{itemize}

It is of interest that extreme points of $f\left( \mathcal{B}\right) $ are
stable under this procedure: 
\begin{equation*}
\mathsf{D}_{2}(\left( l_{2}\right) ^{f})=\left( l_{2}\right) ^{f};\text{ }%
\mathsf{D}_{2}(\left( c_{0}\right) ^{f})=\left( c_{0}\right) ^{f}.
\end{equation*}
\end{remark}

To distinguish between general divisible classes and classes of type $%
\mathsf{D}_{2}(X^{f})$, the last ones will be referred as to \textit{%
2-divisible classes}.

\begin{definition}
A class $X^{f}$ (and its Minkowski's base $\frak{M}(X^{f})$) is said to be
quotient closed if for any $A\in \frak{M}(X^{f})$ and its subspace $%
B\hookrightarrow A$ the quotient $A/B$ belongs to $\frak{M}(X^{f})$.
\end{definition}

Let $K\subseteq \frak{M}$ be a class of finite dimensional Banach spaces
(recall, that isometric spaces are identified). Define operations $H$, $Q$
and $\ast $ that transform a class $K$ in another class - $H(K)$; $Q(K)$ or $%
(K)^{\ast }$ respectively. Namely, let

\begin{align*}
H(K) & =\{A\in\frak{M}:A\hookrightarrow B;\text{ \ }B\in K\} \\
Q(K) & =\{A\in\frak{M}:A=B/F;\text{ \ }F\hookrightarrow B;\text{ \ }B\in K\}
\\
(K)^{\ast} & =\{A^{\ast}\in\frak{M}:A\in K\}
\end{align*}

In words, $H(K)$ consists of all subspaces of spaces from $K$; $Q(K)$
contains all quotient spaces of spaces of $K$; $(K)^{\ast}$ contains all
conjugates of spaces of $K$.

The following theorem lists properties of these operations. In iteration of
the operations parentheses may be omitted.

Thus, $K^{\ast\ast}\overset{def}{=}\left( (K)^{\ast}\right) ^{\ast}$; $HH(K)%
\overset{def}{=}H(H(K))$ and so on.

\begin{theorem}
Any set $K$ of finite dimensional Banach spaces has the following properties:

\begin{enumerate}
\item  $K^{\ast \ast }=K$; $HH(K)=H(K)$; $QQ(K)=Q(K)$;

\item  $K\subset H(K)$; $K\subset Q(K)$;

\item  If $K_{1}\subset K_{2}$ then $H(K_{1})\subset H(K_{2})$ and $%
Q(K_{1})\subset Q(K_{2})$;

\item  $\left( H(K)\right) ^{\ast }=Q(K^{\ast })$; $\left( Q(K)\right)
^{\ast }=H(K^{\ast })$;

\item  $HQ(HQ(K))=HQ(K)$; $QH(QH(K))=QH(K)$.
\end{enumerate}
\end{theorem}

\begin{proof}
1, 2 and 3 are obvious.

4. If $A\in Q(K)$ then $A=B/E$ for some $B\in K$ and its subspace $E$. Thus, 
$A^{\ast }$ is isometric to a subspace of $B^{\ast }.$ Hence, $A^{\ast }\in
H(B^{\ast })$, i.e., $A^{\ast }\in H(K^{\ast })$. Since $A$ is arbitrary, $%
\left( Q(K)\right) ^{\ast }\subseteq H(K^{\ast })$. Analogously, if $B\in K$
and $A\in H(B)$ then $A^{\ast }$ may be identified with a quotient $B^{\ast
}/A^{\perp }$, where $A^{\perp }$ is the annihilator of $A$ in $B^{\ast }$: 
\begin{equation*}
A^{\perp }=\{f\in B^{\ast }:\left( \forall a\in A\right) \text{ \ }\left(
f\left( a\right) =0\right) \}.
\end{equation*}

Hence, $A^{\ast }\in (Q(K^{\ast }))^{\ast }$, and thus $\left( H(K)\right)
^{\ast }\subseteq Q(K^{\ast })$.

From the other hand, 
\begin{align*}
H(K^{\ast })& =\left( H(K^{\ast })\right) ^{\ast \ast }\subseteq \left(
Q(K^{\ast \ast })\right) ^{\ast }=\left( Q(K)\right) ^{\ast }; \\
Q(K^{\ast })& =\left( Q(K^{\ast })\right) ^{\ast \ast }\subseteq \left(
H(K^{\ast \ast })\right) ^{\ast }=\left( H(K)\right) ^{\ast }.
\end{align*}

5. Let $A\in HQ(K)$. Then $A$ is isometric to a subspace of some quotient
space $E/F$, where $E\in K$; $F\hookrightarrow E$. If $B$ is a subspace of $%
A $ then $\left( A/B\right) ^{\ast }=\left( E/F\right) ^{\ast }/B^{\perp }$,
i.e. $\left( Q\left( HQ(K)\right) \right) ^{\ast }\subseteq Q\left(
Q(K)^{\ast }\right) $. Because of 
\begin{align*}
Q\left( HQ(K)\right) & =\left( Q\left( HQ(K)\right) \right) ^{\ast \ast
}\subseteq (Q\left( Q(K)^{\ast }\right) )^{\ast } \\
& \subseteq H(Q(K)^{\ast \ast })=HQ(K),
\end{align*}

we have 
\begin{equation*}
H(Q(H(Q(K))))\subseteq H(H(Q(K)))=HQ(K).
\end{equation*}

Analogously, if $A\in QH(K)$, then $A$ is isometric to a quotient space $F/E$%
, where $F\in H(B)$ for some $B\in K$ and $E\hookrightarrow B$. If $W\in
H(A) $, i.e., if $W\in H(F/E)$, then $W^{\ast }=(F/E)^{\ast }/W^{\perp }$
and $(F/E)^{\ast }$ is isometric to a subspace $E^{\perp }$ of $F^{\ast }\in
(H(B))^{\ast }$. Thus, $(H(QH(K)))^{\ast }\subseteq H((H(K))^{\ast })$ and 
\begin{align*}
H\left( QH(K)\right) & =\left( H\left( QH(K)\right) \right) ^{\ast \ast
}\subseteq (H\left( H(K)^{\ast }\right) )^{\ast } \\
& \subseteq Q(H(K)^{\ast \ast })=QH(K).
\end{align*}

Hence, 
\begin{equation*}
Q(H(Q(HQ(K))))\subseteq Q(Q(H(K)))=QH(K).
\end{equation*}

Converse inclusion follows from 2.
\end{proof}

Consider for a given 2-divisible class $X^{f}$ its Minkowski's base $\frak{M}%
(X^{f})$ and enlarge it by addition to $\frak{M}(X^{f})$ of all quotient
spaces of spaces from $\frak{M}(X^{f})$ and all their subspaces. In the
formal language, consider a set $H(Q(\frak{M}(X^{f})))$.

For any class $W^{f}$ the set $\frak{N}=\frak{M}(W^{f})$ has following
properties:

(\textbf{C}) $\frak{N}\ $is a closed subset of the Minkowski's space $\frak{M%
}$;

(\textbf{H}) If $A\in\frak{N}$ and $B\in H(A)$ then $B\in\frak{N}$;

(\textbf{A}$_{0}$) For any $A$, $B\in\frak{N}$ there exists $C\in\frak{N}$
such that $A\in H(C)$ and $B\in H(C)$.

\begin{theorem}
Let $\frak{N}$ be a set of finite dimensional Banach spaces; $\frak{N}%
\subset \frak{M}$. If $\frak{N}$ has properties (\textbf{C}), (\textbf{H})
and (\textbf{A}$_{0}$) then there exists a class $X^{f}$ such that $\frak{N}=%
\frak{M}(X^{f})$.
\end{theorem}

\begin{proof}
Since (\textbf{A}$_{0}$), all spaces from $\frak{N}$ may be directed in an
inductive isometric system. Its direct limit $W$ generates a class $W^{f}$
with $\frak{M}(W^{f})\supseteq \frak{N}$. As it follows from the property (%
\textbf{C}), $\frak{M}(W^{f})$ contains no spaces besides that of $\frak{N}$%
. Hence, $\frak{M}(W^{f})=\frak{N}$.

Another proof may be given by using an ultraproduct of all spaces from $%
\frak{N}$ by the ultrafilter, coordinated with a partial ordering on $\frak{N%
}$, which is generated by the property (A$_{0}$).
\end{proof}

Let $X$ be a $B$-convex Banach space.

Let $X=l_{2}(Y)$ (and, hence, $X^{f}=\mathsf{D}_{2}(Y^{f})$). Consider the
Minkowski's base $\frak{M}(X^{f})$ and its enlargement $H(Q(\frak{M}%
(X^{f})))=HQ\frak{M}(X^{f})$.

\begin{theorem}
There exists a Banach space $W$ such that $HQ\frak{M}(X^{f})=\frak{M}(W^{f})$%
.
\end{theorem}

\begin{proof}
Obviously, $HQ\frak{M}(X^{f})$ has properties (H) and (C).

Since $\frak{M}(X^{f})$ is 2-divisible, then for any $A,B\in \frak{M}(X^{f})$
a space $A\oplus _{2}B$ belongs to $\frak{M}(X^{f})$ and, hence, to $HQ\frak{%
M}(X^{f})$.

If $A,B\in Q\frak{M}(X^{f})$ then $A=F/F_{1}$; $B=E/E_{1}$ for some $E,F\in Q%
\frak{M}(X^{f})$.

Since $F/F_{1}\oplus _{2}E$ is isometric to a space ($F\oplus
_{2}E)/F_{1}^{\prime }$, where $F_{1}^{\prime }=\{(f,0)\in F\oplus E:f\in
F_{1}\}$, $F/F_{1}\oplus _{2}E$ belongs to $Q\frak{M}(X^{f})$.

Clearly, 
\begin{equation*}
F/F_{1}\oplus _{2}E/E_{1}=(F/F_{1}\oplus _{2}E)/E_{1}^{\prime },
\end{equation*}
where 
\begin{equation*}
E_{1}^{\prime }=\{(o,e)\in F/F_{1}\oplus _{2}E:e\in E_{1}\},
\end{equation*}
and, hence, belongs to $Q\frak{M}(X^{f})$ too.

If $A,B\in HQ\frak{M}(X^{f})$ then $A\hookrightarrow E$, $B\hookrightarrow F$
for some $E,F\in Q\frak{M}(X^{f})$. Since $E\oplus _{2}F\in Q\frak{M}(X^{f})$%
, obviously $A\oplus _{2}B\in HQ\frak{M}(X^{f})$.

Thus, $HQ\frak{M}(X^{f})$ has a property (A$_{0}$). A desired result follows
from the preceding theorem.
\end{proof}

\begin{definition}
Let $X$ be a Banach space, which generates a class $X^{f}$\ of finite
equivalence. A class $\ast \ast (X^{f})$ is defined to be a such class $%
W^{f} $ that $\frak{M}(W^{f})=HQ\frak{M}(\mathsf{D}_{2}Y^{f})$
\end{definition}

Clearly, $W^{f}$ is quotient closed. Obviously, $X^{f}<_{f}W^{f}$. It will
be said that $W^{f}$ is a result of a procedure $\ast \ast $ that acts on $f(%
\mathcal{B})$.

\begin{remark}
The procedure $\ast \ast :X^{f}\rightarrow \ast \ast \left( X^{f}\right)
=W^{f}$ may be regarded as a closure operator on the partially ordered set $%
f\left( \mathcal{B}\right) $. Indeed, it is

\begin{itemize}
\item  Monotone, i.e., $X^{f}<_{f}\ast \ast (X^{f})$;

\item  Idempotent, i.e., $X^{f}=\ast \ast (\ast \ast (X^{f}))$;

\item  Preserve the order: $X^{f}<_{f}Y^{f}$ $\Longrightarrow \ast \ast
(X^{f})<_{f}\ast \ast (Y^{f})$.
\end{itemize}

It is of interest that extreme points of $f\left( \mathcal{B}\right) $ are
stable under this procedure: 
\begin{equation*}
\ast \ast (\left( l_{2}\right) ^{f})=\left( l_{2}\right) ^{f};\text{\ }\ast
\ast (\left( c_{0}\right) ^{f})=\left( c_{0}\right) ^{f}.
\end{equation*}
\end{remark}

\begin{theorem}
For any Banach space $X$ a class $\ast \ast \left( X^{f}\right) $ is
2-divisible.
\end{theorem}

\begin{proof}
Let $\frak{N}=\frak{M}(\ast \ast (X^{f}))$\ . Since for any pair $A,B\in 
\frak{N}$ \ their $l_{2}$-sum belongs to $\frak{N}$, then by an induction, $%
(\sum\nolimits_{i\in I}\oplus A_{i})_{2}\in \frak{N}$ for any finite subset $%
\{A_{i}:i\in I\}\subset \frak{N}$.

Hence, any infinite direct $l_{2}$-sum $(\sum\nolimits_{i\in I}\oplus
A_{i})_{2}$ is finite representable in $\ast \ast \left( X^{f}\right) $. Let 
$\{A_{i}:i<\infty \}\subset \frak{N}$ \ is dense in $\frak{N}$. Let $%
Y_{1}=(\sum\nolimits_{i<\infty }\oplus A_{i})_{2}$; $Y_{n+1}=Y_{n}\oplus
_{2}Y_{1}$; $Y_{\infty }=\underset{\rightarrow }{\lim }Y_{n}$. Clearly, $%
Y_{\infty }=l_{2}\left( Y\right) $ and belongs to $\ast \ast \left(
X^{f}\right) $.
\end{proof}

Let $\star :f\left( \mathcal{B}\right) \rightarrow f\left( \mathcal{B}%
\right) $ be one else procedure that will be given by following steps.

Let $X\in \mathcal{B}$; $Y^{f}=\mathsf{D}_{2}\left( X^{f}\right) $. Let $%
\frak{Y}_{0}$ be a countable dense subset of $\frak{M}(Y^{f})$; $\frak{Y}%
_{0}=\left( Y_{n}\right) _{n<\infty }$. Consider a space $Z=(\sum_{n<\infty
}\oplus Y_{n})_{2}$ and its conjugate $Z^{\ast }$. $Z^{\ast }$ generates a
class $\left( Z^{\ast }\right) ^{f}$ which will be regarded as a result of
acting of a procedure $\star :X^{f}\rightarrow \left( Z^{\ast }\right) ^{f}$%
. Since $\left( Z^{\ast }\right) ^{f}$ is 2-divisible, iterations of the
procedure $\star $ are given by following steps: Let $\frak{Z}_{0}$ is a
countable dense subset of $\frak{M}(\left( Z^{\ast }\right) ^{f})$; $\frak{Z}%
_{0}=\left( Z_{n}\right) _{n<\infty }$. Consider a space $W=(\sum_{n<\infty
}\oplus Z_{n})_{2}$ and its conjugate $W^{\ast }$. Clearly,\ $W^{\ast }$
generates a class 
\begin{equation*}
\left( W^{\ast }\right) ^{f}=\star \left( Z^{\ast }\right) ^{f}=\star \star
\left( X^{f}\right) .
\end{equation*}

\begin{theorem}
For any Banach space $X$ classes $\ast \ast \left( X^{f}\right) $ and $\star
\star \left( X^{f}\right) $ are identical.
\end{theorem}

\begin{proof}
From the construction follows that $H(Z^{\ast })=\left( QH(l_{2}(X)\right)
^{\ast }$ and that 
\begin{align*}
H(W^{\ast })& =(QH(Z^{\ast }))^{\ast }=(Q(QH(l_{2}(X)))^{\ast })^{\ast } \\
& =H(QH(l_{2}(X)))^{\ast \ast }=HQH(l_{2}(X)).
\end{align*}

Hence, 
\begin{equation*}
\frak{M}(\left( W^{\ast }\right) ^{f})=HQ\frak{M}(Y^{f})=HQ\frak{M}(\mathsf{D%
}_{2}(X^{f})).
\end{equation*}
\end{proof}

\begin{theorem}
Let $X$ be a Banach space, which generates a class of finite equivalence $%
X^{f}$. If $X$ is B-convex, then the procedure $\ast \ast $ maps $X^{f}$ to
a class $\ast \ast \left( X^{f}\right) $ with the same lower and upper
bounds of its $l_{p}$-spectrum as $X^{f}$. If $X$ is superreflexive then $%
\ast \ast \left( X^{f}\right) $ is superreflexive too.
\end{theorem}

\begin{proof}
Obviously, if $p\in S(X)$ and $p\leq 2$ then the whole interval $[p,2]$
(that may consist of one point) belongs to $S(X)$ and, hence, to $S(\ast
\ast \left( X^{f}\right) )$. If $p\in S(X)$ and $p>2$ then, by duality, $%
p/(p-1)\in S\left( \star \left( X^{f}\right) \right) $; hence $%
[p/(p-1),2]\subset S\left( \star \left( X^{f}\right) \right) $ and, thus $%
[2,p]\subset S(X)\subset S(\ast \ast \left( X^{f}\right) )$. If $p\notin
S(X) $ then $p/(p-1)\notin S\left( \star \left( X^{f}\right) \right) $ by
its construction and, hence, $p\notin S(\ast \ast \left( X^{f}\right) )$ by
the same reason. Hence, if 
\begin{equation*}
p_{X}=\inf S(X)\text{; \ }q_{X}=\sup S(X)
\end{equation*}
then 
\begin{equation*}
S(\ast \ast \left( X^{f}\right) )=[p_{X},q_{X}]=[\inf S(X),\sup S(X)].
\end{equation*}

The second assertion of the theorem is also obvious.
\end{proof}

\begin{remark}
If $X$ is not $B$-convex, then 
\begin{equation*}
\star \left( X^{f}\right) =\star \star \left( X^{f}\right) =\ast \ast \left(
X^{f}\right) =\left( c_{0}\right) ^{f}.
\end{equation*}
\end{remark}

\section{Spaces of almost universal disposition}

Recall a definition, which is due to V.I. Gurarii [7].

\begin{definition}
Let $X$ be a Banach space; $\mathcal{K}$ be a class of Banach spaces. $X$ is
said to be a space of almost universal disposition with respect to $\mathcal{%
K}$ if for any pair of spaces $A$, $B$ of $\mathcal{K}$ such that $A$ is a
subspace of $B$ ($A\hookrightarrow B$), every $\varepsilon >0$ and every
isomorphic embedding $i:A\rightarrow X$ there exists an isomorphic embedding 
$\hat{\imath}:B\rightarrow X$, which extends $i$ (i.e., $\hat{\imath}|_{A}=i$%
) and such, then 
\begin{equation*}
\left\| \hat{\imath}\right\| \left\| \hat{\imath}^{-1}\right\| \leq
(1+\varepsilon )\left\| i\right\| \left\| i^{-1}\right\| .
\end{equation*}
\end{definition}

In the construction of the classical Gurarii's space of almost universal
disposition with respect to $\frak{M}$ \ the main role plays a property of $%
\frak{M}$, which may be called the \textit{isomorphic amalgamation property}.

\begin{definition}
Let $X\in \mathcal{B}$ \ generates a class $X^{f}$\ with a Minkowski's base $%
\frak{M}(X^{f})$. It will be said that $\frak{M}(X^{f})$ (and the whole
class $X^{f}$) has the isomorphic amalgamation property if for any fifth $%
\left\langle A,B_{1},B_{2},i_{1},i_{2}\right\rangle $, where $A$, $B_{1}$, $%
B_{2}\in \frak{M}(X^{f})$; $i_{1}:A\rightarrow B_{1}$ and$\
i_{2}:A\rightarrow B_{2}$ are isomorphic embeddings, there exists a triple $%
\left\langle j_{1},j_{2},F\right\rangle $, where $F\in \frak{M}(X^{f})$ and$%
\ j_{1}:B_{1}\rightarrow F$; $j_{2}:B_{2}\rightarrow F$ are isometric
embedding, such that $j_{1}\circ i_{1}=j_{2}\circ i_{2}$.
\end{definition}

\begin{theorem}
Any 2-divisible quotient closed class $X^{f}$ has the isomorphic
amalgamation property.
\end{theorem}

\begin{proof}
Let $A$, $B_{1}$, $B_{2}\in \frak{M}(X^{f})$; $i_{1}:A\rightarrow B_{1}$ and 
$i_{2}:A\rightarrow B_{2}$ are isomorphic embeddings. Since $X^{f}$ is
2-divisible, $C=B_{1}\oplus _{2}B_{2}\in \frak{M}(X^{f})$. Consider a
subspace $H$ of $C$ that is formed by elements of kind $\left(
i_{1}a,-i_{2}a\right) $, where $a$ runs $A$: 
\begin{equation*}
H=\{\left( i_{1}a,-i_{2}a\right) :a\in A\}.
\end{equation*}

Consider a quotient $W=C/H$. Since $X^{f}$ is quotient closed, $W\in \frak{M}%
(X^{f})$. Let $h:C\rightarrow W$ be a standard quotient map. Let $%
j_{1}:B_{1}\rightarrow W$ and $j_{2}:B_{2}\rightarrow W$ are given by: 
\begin{align*}
j_{1}(b_{1})& =h(b_{1},0)\text{ \ \ \ for }b_{1}\in B_{1}; \\
j_{2}(b_{2})& =h(0,b_{2})\text{ \ \ \ for }b_{2}\in B_{2}.
\end{align*}

It is clear that $j_{1}$ and $j_{2}$ are isometric embeddings such that $\
j_{1}\circ i_{1}=j_{2}\circ i_{2}$.
\end{proof}

\begin{remark}
It may be shown that any quotient closed divisible (not necessary
2-divisible) class $X^{f}$ enjoys the isomorphic amalgamation property too.
For aims of the article the preceding result is satisfactory.
\end{remark}

The proof of a following results is almost literally repeats the Gurarii's
proof [7] on existence of a space of almost universal disposition with
respect to $\frak{M}$. Only changes that need to be made are: a substitution
of a set $\frak{M}$ with $\frak{M}(X^{f})$ for a given class $X^{f}$ and
using instead the mentioned above isomorphic amalgamation property of a set $%
\frak{M}$ \ the same property of a set $\frak{M}(X^{f})$. For this reason
the proof of it is omitted.

\begin{theorem}
Any $\ast \ast $-closed class $X^{f}$ contains a separable space $E_{X}$ of
almost universal disposition with respect to a set $\frak{M}(X^{f})$. This
space is unique up to almost isometry and is almost isotropic (in an
equivalent terminology, has an almost transitive norm: for any two elements $%
a$, $b\in E_{X}$, such that $\left\| a\right\| =\left\| b\right\| $ \ and
every $\varepsilon >0$ there exists an automorphism $u=u(a,b,\varepsilon
):E_{X}\overset{onto}{\rightarrow }E_{X}$ such that $\left\| u\right\|
\left\| u^{-1}\right\| \leq 1+\varepsilon $ and $ua=b$). This space is an
approximative envelope of a class $X^{f}$: for any $\varepsilon >0$ every
separable Banach space which is finitely representable in $X^{f}$ is $%
(1+\varepsilon )$-isomorphic to a subspace of $E_{X}$.
\end{theorem}

\section{Subspace Homogeneous and Quotient Homogeneous Banach Spaces}

\begin{theorem}
Let $X^{f}=\ast \ast \left( X^{f}\right) $ be a superreflexive class; $E_{X}$
be a corresponding separable space of almost universal disposition. Then $%
E_{X}$ is subspace homogeneous.
\end{theorem}

\begin{proof}
Let $Z_{1}$ and $Z_{2}$ be subspaces of $E_{X}$ that are of equal
codimension in $E_{X}$ (i.e. $codim\left( Z_{1}\right)
=dim(E_{X}/Z_{1})=dim(E_{X}/Z_{2})=codim\left( Z_{2}\right) $).

If $codim\left( Z_{1}\right) <\omega $ (here and below $\omega $ denotes the
first infinite ordinal) then assertions of the theorem are satisfied
independently of a property of $E_{X}$ to be a space of almost universal
disposition. Indeed, for an arbitrary Banach space $X$ any its subspaces of
equal codimension, say, $Y$ and $Z$, are isomorphic and each of them has a
complement in $X$. So, $X$ may be represented either as a direct sum $%
X=Y\oplus A$ or as $X=Z\oplus B$ where $dim(A)=dim(B)<\omega $.

Let $v:A\rightarrow B$ and $u:Y\rightarrow Z$ be isomorphisms. Clearly, $%
v\oplus u:X\rightarrow X$ be an automorphism of $X$ that extends $u$.

Assume that $Z_{1}$ and $Z_{2}$ are isomorphic subspaces of $E_{X}$ of
infinite codimension. Let $u:Z_{1}\rightarrow Z_{2}$ be the corresponding
isomorphism.

Let $Z_{1}$ be represented as a direct limit (= as a closure of the union)
of a chain 
\begin{equation*}
Z_{1}^{\prime }\hookrightarrow Z_{2}^{\prime }\hookrightarrow
...\hookrightarrow Z_{n}^{\prime }\hookrightarrow ...\hookrightarrow Z_{1};
\end{equation*}
namely, $Z_{1}=\overline{\cup \{Z_{n}^{\prime }:n<\omega \}}$ and $\cup
Z_{n}^{\prime }$ is dense in $Z_{1}$. Clearly, 
\begin{equation*}
uZ_{1}^{\prime }\hookrightarrow uZ_{2}^{\prime }\hookrightarrow
...\hookrightarrow uZ_{n}^{\prime }\hookrightarrow ...\hookrightarrow
uZ_{1}=Z_{2}
\end{equation*}
and the union $\cup uZ_{n}^{\prime }$ is dense in $Z_{2}$.

Let $\left( e_{n}\right) _{n<\omega }$ be a sequence of linearly independent
elements of $E_{X}$ of norm one, which linear span is dense in $E_{X}$.

Denote a restriction $u\mid _{Z_{n}^{\prime }}$ by $u_{n}$and define two
sequences $\left( f_{n}\right) $ and $\left( g_{n}\right) $ of elements of $%
E_{X}$ and a sequence $\left( v_{n}\right) $ of isomorphisms by an induction.

Let $f_{1}$ be an element of $\left( e_{n}\right) _{n<\omega }$ with a
minimal number, which does not belongs to $Z_{1}^{\prime }$. Let $%
W_{1}=span\left( Z_{1}^{\prime }\cup \{f_{1}\}\right) $ (as $span(A)$ here
and below is denoted a closure of the linear span of $A$). Let $\varepsilon
>0$ and $v_{1}:W_{1}\rightarrow E_{X}$ be an extension of $u_{1}$ (i.e. $%
v_{1}\mid _{W_{1}}=u_{1}$) such that $\left\| v_{1}\right\| \left\|
v_{1}^{-1}\right\| \leq \left( 1+\varepsilon \right) \left\| u_{1}\right\|
\left\| u_{1}^{-1}\right\| $. Let $g_{1}=v_{1}f_{1}$; $U_{1}=v_{1}W_{1}$.

Let $g_{2}$ be an element of $\left( e_{n}\right) _{n<\omega }$ with a
minimal number, which does not belongs to $span\left( uZ_{1}^{\prime }\cup
U_{1}\right) $ and $U_{2}=span\left( uZ_{1}^{\prime }\cup U_{1}\cup
\{g_{2}\}\right) $. Since $E_{X}$ is a space of almost universal
disposition, there exists an isomorphism $\left( v_{2}\right) ^{-1}$ which
extends both $\left( u_{2}\right) ^{-1}$ and $\left( v_{1}\right) ^{-1}$ and
satisfies 
\begin{equation*}
\left\| v_{2}\right\| \left\| v_{2}^{-1}\right\| \leq \left( 1+\varepsilon
^{2}\right) \max \{\left\| u_{2}\right\| \left\| u_{2}^{-1}\right\| ;\left\|
v_{1}\right\| \left\| v_{1}^{-1}\right\| \}.
\end{equation*}
Let $f_{2}=\left( v_{2}\right) ^{-1}g_{2}$. This close the first step of the
induction.

Assume that $\{f_{1},f_{2},...,f_{n}\}$; $\{g_{1},g_{2},...,g_{n}\}$; $%
\{v_{1},v_{2},...,v_{n}\}$; $\{U_{1},U_{2},...,U_{n}\}$ and $%
\{W_{1},W_{2},...,W_{n}\}$ are already chosen. If $n$ is odd, we choose
sequentially:

$f_{n+1}$as an element of $\left( e_{n}\right) _{n<\omega }$ with a minimal
number, which does not belongs to $Z_{n+1}^{\prime }$; $W_{n+1}=span\left(
Z_{n+1}^{\prime }\cup \{f_{n+1}\}\right) $; $v_{n+1}:W_{n+1}\rightarrow
E_{X} $ be an extension of $u_{n+1}$ such that 
\begin{equation*}
\left\| v_{n+1}\right\| \left\| v_{n+1}^{-1}\right\| \leq \left(
1+\varepsilon ^{n+1}\right) \left\| u_{n+1}\right\| \left\|
u_{n+1}^{-1}\right\| ;
\end{equation*}
$g_{n+1}=v_{n+1}f_{n+1}$; $U_{n+1}=v_{n+1}W_{n+1}$.

If $n$ is even, we choose sequentially:

$g_{n+1}$ as an element of $\left( e_{n}\right) _{n<\omega }$ with a minimal
number, which does not belongs to $span\left( uZ_{n}^{\prime }\cup
U_{n}\right) $; $U_{n+1}=span\left( uZ_{n}^{\prime }\cup U_{n}\cup
\{g_{n+1}\}\right) $; $\left( v_{n+1}\right) ^{-1}$which extends both $%
\left( u_{n+1}\right) ^{-1}$ and $\left( v_{n}\right) ^{-1}$ and satisfies 
\begin{equation*}
\left\| v_{n+1}\right\| \left\| v_{n+1}^{-1}\right\| \leq \left(
1+\varepsilon ^{2}\right) \max \{\left\| u_{n+1}\right\| \left\|
u_{n+1}^{-1}\right\| ;\left\| v_{n}\right\| \left\| v_{n}^{-1}\right\| \};
\end{equation*}
$f_{n+1}=\left( v_{n+1}\right) ^{-1}g_{n+1}$

Since $E_{X}$ is superreflexive, a sequence of isomorphisms $\left(
v_{n}\right) $ converges to an automorphism $V:E_{X}\rightarrow E_{X}$,
which satisfies 
\begin{equation*}
\left\| V\right\| \left\| V^{-1}\right\| \leq \prod\nolimits_{n=1}^{\infty
}\left( 1+\varepsilon ^{n}\right) \left\| u\right\| \left\| u^{-1}\right\| .
\end{equation*}
\end{proof}

\begin{corollary}
Let $X^{f}=\ast \ast \left( X^{f}\right) $ be a superreflexive class; $E_{X}$
be a corresponding separable space of almost universal disposition. Then $%
\left( E_{X}\right) ^{\ast }$ is quotient homogeneous.
\end{corollary}

\begin{proof}
The desired result is a consequence of the superreflexivity and follows by
duality.
\end{proof}

\section{The Banach problem on squares}

A classical Banach's problem on squares (whether every Banach space is
isomorphic to its Carthesian square?) was solved in [2]. At the same time
there appeared another solution [3]; one more counterexample was presented
in [8].

Superreflexive spaces of almost universal disposition form a wide class of
spaces that are not isomorphic to their Carthesian squares, which is
different from listed above.

\begin{theorem}
Let $\ast \ast X^{f}$ be a superreflexive class; $E_{X}$ be a corresponding
separable space of almost universal disposition. Then either $E_{X}$ is
isomorphic to the Hilbert space $l_{2}$ or spaces $E_{X}$ and $E_{X}\oplus
E_{X}$ are not isomorphic.
\end{theorem}

\begin{proof}
Let the symbol $Y\approx Z$ means that $Y$ and $Z$ are isomorphic.

Assume that $E_{X}\approx E_{X}\oplus E_{X}$. Then $E_{X}^{\ast }\approx
E_{X}^{\ast }\oplus E_{X}^{\ast }$. Since $E_{X}^{\ast }$ is not isomorphic
to $l_{2}$, there exists a subspace $W^{\ast }$ of $E_{X}^{\ast }$ such that 
$E_{X}^{\ast }\oplus W^{\ast }$ and $W^{\ast }$ are not isomorphic. Clearly, 
$W^{\ast }$ may be identified with $\left( E_{X}/Y\right) ^{\ast }$ for some 
$Y\hookrightarrow E_{X}$. Certainly, $E_{X}\oplus \left( E_{X}/Y\right) $
and $E_{X}/Y$ are not isomorphic.

Let $u:E_{X}\rightarrow E_{X}\oplus E_{X}$ be an isomorphism. Then $%
uY\hookrightarrow E_{X}\oplus E_{X}$ and $uY\approx Y$.

Consider two isomorphic embeddings, say, $j_{1}$ and $j_{2}$ of $Y$ into $%
E_{X}$.

Let $j_{1}:Y\hookrightarrow E_{X}$ be the identical embedding. Let $%
j_{2}=u\circ j_{1}$.

Let $v:j_{1}Y\rightarrow j_{2}Y$ be an isomorphism. Then $\left( v^{-1}\circ
j_{2}\right) Y\hookrightarrow E_{X}$. However the isomorphism $v$ cannot be
extended to any automorphism of $E$ because of absence of isomorphisms
between $E_{X}\oplus \left( E_{X}/Y\right) $ and $E_{X}/Y$. But this
contradicts with the property of $E_{X}$ to be subspace homogeneous.
\end{proof}

\begin{remark}
If a Banach space $X$ is not isomorphic to any space $Y$ that generates a
divisible class $Y^{f}$, then any space $Z\in X^{f}$ is not isomorphic to
its Carthesian square. First example of such $X$ was the space $F=\left(
\sum \oplus l_{p_{i}}^{\left( n_{i}\right) }\right) _{2}$, where $%
p_{i}\searrow 2$; $n_{i}=dim\left( l_{p_{i}}^{\left( n_{i}\right) }\right)
\rightarrow \infty $; $\left( p_{i}-2\right) \log n_{i}\rightarrow \infty $,
that was constructed by T. Figiel [8].
\end{remark}

\section{The Schroeder-Bernstein problem}

The \textit{Schroeder-Bernstein problem,} that was induced by the classical
paper [9] and received its name thanks to a certain analogy may be
formulated as follows:

\textit{Let }$X$\textit{, }$Y$\textit{, }$Z$\textit{\ be Banach spaces. Is
there true that a condition }$X\approx X\oplus Y\oplus Z$\textit{\ implies
that }$X\approx X\oplus Y$\textit{?}

Recall that in [9] was shown that if $X\approx X\oplus X$ and $Y\approx
Y\oplus Y$ then $X\approx X\oplus Y\oplus Z$ implies that $X\approx X\oplus
Y $.

This problem was (in negative) solved by W.T. Gowers (see [5]).

An another solution of this problem based on the property of spaces of kind $%
E_{X}$, which needs additional definitions.

\begin{definition}
(Cf.[10]). Let $X$ be a Banach space; $Y$ - its subspace. $Y$ is said to be
a reflecting subspace of $X$ (symbolically: $Y\prec _{u}X$) if for every $%
\varepsilon >0$ and every finite dimensional subspace $A\hookrightarrow X$
there exist an isomorphic embedding $u:A\rightarrow Y$ such that $\left\|
u\right\| \left\| u^{-1}\right\| \leq 1+\varepsilon $ and $u\mid _{A\cap
Y}=Id_{A\cap Y}$.
\end{definition}

As was shown in [10], if $Y\prec _{u}X$ then $Y^{\ast \ast }$ is an image of
a norm one projection $P:X^{\ast \ast }\rightarrow Y^{\ast \ast }$ (under
canonical embedding of $Y^{\ast \ast }\ $into $X^{\ast \ast }$).

\begin{definition}
(Cf\textrm{\ }[11]). A Banach space $E$ is said to be existentialy closed in
a class $X^{f}$ if for any isometric embedding $i:E\rightarrow Z$ into an
arbitrary space $Z\in X^{f}$ its image $iE$ is a reflecting subspace of $Z$: 
$iY\prec _{u}Z$.
\end{definition}

A class of all spaces $E$ that are existentialy closed in $X^{f}$ is denoted
by $\mathcal{E}\left( X^{f}\right) $. In [11] was shown that for any Banach
space $X$ the class $\mathcal{E}\left( X^{f}\right) $ is nonempty; moreover,
any $Y<_{f}X^{f}$ may be isometricaly embedded into some $E\in \mathcal{E}%
\left( X^{f}\right) $ of the dimension $dim(E)=\max \{dim(Y),\omega \}$.

\begin{theorem}
For any class $\ast \ast (X^{f})$ the corresponding space $E_{X}$ of almost
universal disposition is existentialy closed in $\ast \ast (X^{f})$.
\end{theorem}

\begin{proof}
Let $E_{X}\hookrightarrow Z\in \ast \ast (X^{f})$. Let $A\hookrightarrow Z$.
Denote $E_{X}\cap A=B$. Since $E_{X}\sim _{f}Z$, for any $\varepsilon >0$
there exists an isomorphic embedding $u:A\rightarrow E_{X}$ such that $%
\left\| u\right\| \left\| u^{-1}\right\| \leq 1+\varepsilon /4$. Let $F=uA$; 
$G=u\left( E_{X}\cap A\right) =uB$. Since $F\subset E_{X}$ and $G\subset
E_{X}$, a restriction $u\mid _{B}=v$ may be extended to an isomorphism $%
V:F\rightarrow E_{X}$ such that $\left\| V\right\| \left\| V^{-1}\right\|
\leq 1+\varepsilon /4$. Clearly, $VF$ is $(1+\varepsilon )$-isomorphic to $A$
and $VF\cap A=B$. Because of $A$ is arbitrary, $E_{X}\prec _{u}Z$ and,
hence, $E_{X}\in \mathcal{E}\left( \ast \ast (X^{f})\right) $.
\end{proof}

\begin{corollary}
If $\ast \ast \left( X^{f}\right) $ is superreflexive, then $E_{X}$ is a
norm one complemented subspace of any space $Z$ from the class $\ast \ast
\left( X^{f}\right) $ that it contain.
\end{corollary}

\begin{proof}
If $E_{X}\hookrightarrow Z\in \ast \ast (X^{f})$ then $E_{X}\prec _{u}Z$ and
hence $\left( E_{X}\right) ^{\ast \ast }=E_{X}$ is a norm one complemented
subspace of $Z^{\ast \ast }=Z$ according to [10].
\end{proof}

\begin{corollary}
Thou $E_{X}\approx E_{X}\oplus E_{X}\oplus Z$ but $E_{X}\ $and $E_{X}\oplus
E_{X}$ are not isomorphic.
\end{corollary}

\begin{proof}
Because $E_{X}\oplus _{2}E_{X}$ is finitely equivalent to $E_{X}$, it may be
isomorphicaly embedded into $E_{X}$. Since the image of $E_{X}$ in $E_{X}$
is a complemented subspace, clearly, $E_{X}\oplus E_{X}$ is isomorphic to a
complemented subspace of $E_{X}$, i.e., $E_{X}\approx E_{X}\oplus
E_{X}\oplus Z$. In a previous section was proved that $E_{X}\ $and $%
E_{X}\oplus E_{X}$ are not isomorphic.
\end{proof}

\section{The approximation problem}

\begin{definition}
A Banach space $X$ has the \textit{approximation property}\textbf{\ (}shortly%
\textbf{: }AP) provided for every compact subset $K\subset X$ and every $%
\varepsilon >0$ there exists a (bounded linear) operator $u_{K}:X\rightarrow
X$ with a finite dimensional range such that $\left\| u_{K}x-x\right\| \leq
\varepsilon $ for all $x\in K$.

If there exists a constant $\lambda $, independent on $K$ and $\varepsilon $%
, such that norms of operators $u_{K}$ are uniformly bounded by $\lambda $: $%
\left\| u_{K}\right\| \leq \lambda $, then $X$ is said to has the $\lambda $-%
\textit{bounded approximation property}\textbf{\ }($\lambda $-BAP$)$.

If $X$ has the $\lambda $-BAP for all $\lambda >1,$ then $X$ is said to has
the \textit{metric approximation property}\textbf{\ }(MAP).
\end{definition}

The approximation problem, i.e. a question: ''whether every Banach space has
the approximation property?'' was solved (in negative) by P. Enflo [2]. This
solution ( and other known counterexamples to the approximation problem) are
quite complicated.

This section is devoted to prove that a lot of spaces of kind $E_{X}$ - i.e.
separable spaces of almost universal disposition, that are presented in
every quotient closed and divisible class $X^{f}$ - does not enjoy the
approximation property. Unfortunately, using methods cannot give a chance to
prove the absence of approximation property in a non reflexive case. It
seems that result presented below also is far from finality.

\begin{theorem}
Let $\ast \ast \left( X^{f}\right) $ be a superreflexive class such that its 
$l_{p}$-spectrum $S\left( \ast \ast \left( X^{f}\right) \right) $ contains a
point $p\neq 2$. Let $E_{X}$ be the corresponding space of almost universal
disposition. Then $E_{X}$ does not have the approximation property.
\end{theorem}

\begin{proof}
Let $p\in S\left( \ast \ast \left( X^{f}\right) \right) $ and $p\neq 2$.
Then either $p_{X}=\inf S\left( \ast \ast \left( X^{f}\right) \right) <2$ or 
$q_{X}=\sup S\left( \ast \ast \left( X^{f}\right) \right) >2$. From [12]
follows that if $r=p_{X}$ or $r=q_{X}$ then $E_{X}$ (as any other space from 
$\ast \ast \left( X^{f}\right) $) contains a sequence $\left( Y_{n}\right)
_{n<\infty }$ of finite dimensional subspaces of unboundedly growing
dimension that are uniformly complemented in $E_{X}$ and are uniformly
isomorphic to $l_{r}^{k}$. More precise, $dim(Y_{n})\rightarrow \infty $ and
there exists such a constant $C$ that $\sup d(Y_{n},l_{r}^{dim(Y_{n})})<C$
and that 
\begin{equation*}
\sup \{\left\| P_{n}\right\| :P_{n}:E_{X}\rightarrow Y_{n}\}
\end{equation*}
for some sequence of projections $\left( P_{n}\right) $.

Since $E_{X}$ is a\ space of almost universal disposition, it may be assumed
without loss of generality that the sequence $\left( Y_{n}\right) _{n<\infty
}$ is directed in a chain 
\begin{equation*}
Y_{1}\hookrightarrow Y_{2}\hookrightarrow ...\hookrightarrow
Y_{n}\hookrightarrow ...\hookrightarrow E_{X}.
\end{equation*}

From [13] follows that $Y=\overline{\cup _{n}Y_{n}}$ is an $\mathcal{L}_{p}$%
-space which cannot has the complement in $E_{X}$ if $r\neq 2$. Indeed, if $Y
$ is a complemented subspace of $E_{X}$ then $E_{X}$ contains a complemented
subspace isomorphic to $l_{r}$. At the same time, $E_{X}$ contains an
uncomplemented subspace, isomorphic to $l_{r}$. Indeed (we assume that $%
r\neq 2$), as an universal space for a class of separable spaces from $%
\left( \ast \ast \left( X^{f}\right) \right) ^{<f}$, $E_{X}$ contains a
subspace isomorphic to $L_{r}[0,1]$ (for simplicity assume that $L_{r}[0,1]$
is a subspace of $E_{X}$ itself). The last one contains an uncomplemented
subspace $Z$ isomorphic to $l_{r}$ (for $p>2$ this fact is contained in
[14]; for $p<2$ follows from [15]). If there is a projection from $E_{X}$
onto $Z$ then its restriction is a projection from $L_{r}[0,1]$ onto $Z$
what is impossible. Thus, $E_{X}$ contains an uncomplemented subspace,
isomorphic to $l_{r}$. Since $E_{X}$ is subspace homogeneous, every its
subspace, isomorphic to $l_{r}$ has no complement in $E_{X}$.

Assume that $E_{X}$ has the approximation property. By the superreflexivity, 
$E_{X}$ has the metric approximation property. Hence, there exists a
sequence $\left( E^{\left( n\right) }\right) $ of finite dimensional
subspaces of $E_{X}$ and a sequence of operators $R_{n}:E_{X}\rightarrow
E^{\left( n\right) }$ with $\sup \left\| R_{n}\right\| \leq \lambda $ for a
given $\lambda >1$, which convergents pointwise to the identical operator $%
Id_{E_{X}}$. It may be assumed that $Y_{n}\subset E^{\left( m_{n}\right) }$
for a some sequence $\left( m_{n}\right) $ of natural numbers. A restriction 
$P_{n}^{\prime }$ of$\ P_{n}$ to $E^{\left( m_{n}\right) }$ is a projection
from $E^{\left( m_{n}\right) }$ onto $Y_{n}$. A sequence $\left(
R_{m_{m}}\circ P_{n}^{\prime }\right) $ convergents pointwise to a
projection $P:E_{X}\rightarrow Y=\overline{\cup _{n}Y_{n}}$. As it was
shown, this is impossible.
\end{proof}

\begin{problem}
Whether non reflexive space of kind $E_{X}$ either is isomorphic
(=isometric) to the classical Gurarii space of almost universal disposition
or does not enjoy the approximation property?
\end{problem}

\section{Tsirelson like spaces}

A Banach space $X$ is said to has the Tsirelson property if it contains no
subspaces isomorphic to the spaces $l_{p}$ $(1\leq p<\infty )$ and $c_{0}$.

The first example of a Banach space with a such property was constructed by
B.S. Tsirelson [16]. Below it will be shown that spaces that enjoy the
Tsirelson property may be found among spaces of kind $E_{X}$.

\begin{theorem}
Let $\ast \ast \left( X^{f}\right) $ be a superreflexive class such that its 
$l_{p}$-spectrum $S\left( \ast \ast \left( X^{f}\right) \right) $ is a
single point $p=2$: $S\left( \ast \ast \left( X^{f}\right) \right) =\{2\}$.
Let $E_{X}$ be the corresponding space of almost universal disposition. If $%
E_{X}$ is not of cotype 2 then either its conjugate $\left( E_{X}\right)
^{\ast }$ has the Tsirelson property or $E_{X}$ is isomorphic to the Hilbert
space.
\end{theorem}

\begin{proof}
Assume that $E_{X}$ is not isomorphic to the Hilbert space.

Since $E_{X}$ is not of cotype 2, according to the Kwapie\'{n}'s result [17]
its 2-type constant $T_{2}\left( X\right) $ is infinite. So, for any number $%
N$ there exists a finite dimensional subspace $X^{\left( n\right) }$ with $%
T_{2}\left( X^{\left( n\right) }\right) $, where, certainly, $%
n=dim(X^{\left( n\right) })=n(N)$ depends on $N$ (and tends to infinity with 
$N$).

Fix $\varepsilon >0$; $1<C<\infty $ and, proceeding by induction, chose
sequences $\left( p_{k}\right) \subset \mathbb{R}_{+}$; $p_{k}<p_{k+1}<2$; $%
\left( n_{k}\right) \subset \mathbb{N}$ and $\left( N_{k}\right) \subset 
\mathbb{R}$.

Let $N_{1}>2$ be arbitrary chosen. Let $n_{1}$ be the least natural number
such that there exists a subspace $X^{\left( n_{1}\right) }\hookrightarrow X$
of dimension $n_{1}$ for which 
\begin{equation*}
T_{2}(X^{\left( n_{1}\right) })\geq N_{1}.
\end{equation*}

Let 
\begin{equation*}
p_{1}=2-C\left( \log \left( n_{1}\right) \right) ^{-1}.
\end{equation*}

As $N_{2}$ may be chosen any number such that 
\begin{equation*}
\left( N_{2}\right) ^{C\left( \log \left( n_{1}\right) \right) ^{-1}}\geq
2N_{1}.
\end{equation*}

Assume that $N_{1}$, $N_{2}$, ..., $N_{k}$; $p_{1}$, $p_{2}$, ..., $p_{k-1}$%
; $n_{1}$, $n_{2}$, ..., $n_{k-1}$ and subspaces $X^{\left( n_{1}\right) }$, 
$X^{\left( n_{2}\right) }$, ..., $X^{\left( n_{k-1}\right) }$ are chosen.
Let $n_{k}$ be the least natural number such that there exists a subspace $%
X^{\left( n_{k}\right) }\hookrightarrow X$ of dimension $n_{k}$ for which 
\begin{equation*}
T_{2}(X^{\left( n_{k}\right) })\geq N_{k}.
\end{equation*}

Put 
\begin{equation*}
p_{k}=2-C\left( \log \left( n_{k}\right) \right) ^{-1}.
\end{equation*}

Let $N_{k+1}$ satisfies 
\begin{equation*}
\left( N_{k+1}\right) ^{C\left( \log \left( n_{k}\right) \right) ^{-1}}\geq
2N_{k}.
\end{equation*}

It was noted before that $L_{p}$ for $p<2$ has an uncomplemented subspace,
isomorphic to $l_{2}$. The finite-dimensional version of this fact is:

\textit{There exists such a constant }$c$\textit{\ that for every }$N>1$%
\textit{\ there exists such }$n=n(N)$\textit{\ that }$l_{p}^{\left( n\right)
}$\textit{\ contains a subspace }$Y^{\left( m_{n}\right) }$\textit{\ of
dimension }$m_{n}$\textit{\ such that } 
\begin{equation*}
d(Y^{\left( m_{n}\right) },l_{2}^{\left( m_{n}\right) })\leq c\text{; \ }%
\left\| P\right\| \geq N\text{\textit{for every projection} }P:l_{p}^{\left(
n\right) }\rightarrow Y^{\left( m_{n}\right) }.
\end{equation*}

From [18] follows that for any $p\in \left( 1,2\right) $ every Banach space $%
X$, which type $p$-constant is $T_{p}(X)$, contains a subspace, which is $%
(1+\varepsilon )$-isomorphic to $l_{p}^{\left( m\right) }$, of dimension $%
m=m(p,\varepsilon )$ that is estimated as 
\begin{equation*}
m=m(p,\varepsilon )\geq C(\varepsilon )\left( T_{p}(X)\right) ^{1/q},
\end{equation*}
where $p^{-1}+q^{-1}=1$ and $C(\varepsilon )$ depends only on $\varepsilon $.

From the definition of type $p$-constants follows that for a finite
dimensional Banach space $X$ 
\begin{equation*}
\left| T_{p}\left( X\right) -T_{2}\left( X\right) \right| \leq \left(
dim(X)\right) ^{1/p-1/2}.
\end{equation*}

At least, according to the Kwapie\'{n}'s result, for any Banach space $X$ of
cotype 2 
\begin{equation*}
T_{p}(X)\asymp d(X,l_{2}^{dim\left( X\right) }).
\end{equation*}

Since our choosing of $\left\langle N_{k},n_{k},p_{k},X^{\left( n_{k}\right)
}\right\rangle $'s follows:

\begin{enumerate}
\item  $\left| T_{p}\left( X^{\left( n_{k}\right) }\right) -T_{2}\left(
X^{\left( n_{k}\right) }\right) \right| \leq C$ for all $k<\infty $;

\item  $X^{\left( n_{k}\right) }$ contains a subspace that is $%
(1+\varepsilon )$-isomorphic to $l_{p}^{\left( m(p_{k},\varepsilon )\right) }
$.
\end{enumerate}

Let $m_{k}=m(p_{k},\varepsilon )$. Then 
\begin{eqnarray*}
d\left( l_{p_{k}}^{n_{k}},l_{2}^{n_{k}}\right)  &=&\left( m_{k}\right)
^{1/p_{k}-1/2}\geq \left( C(\varepsilon )\left( T_{p_{k}}(X^{\left(
n_{k}\right) })\right) ^{1/q}\right) ^{1/p_{k}-1/2} \\
&\geq &C(\varepsilon )\left( T_{2}(X^{\left( n_{k}\right) })-C\right)
^{1/q}\geq C(\varepsilon )\left( 2^{k}N_{1}-C\right) ^{1/q}
\end{eqnarray*}
and, hence, tends to infinity with $k$.

So, $E_{X}$ (as any other Banach space $X$ of cotype 2 and the single-point
spectrum $S(X)=\{2\}$) contains a sequence of finite dimensional subspaces $%
\left( Y^{m_{k}}\right) $, which are uniformly isomorphic to Euclidean ones
and whose projection constants 
\begin{equation*}
\lambda \left( Y^{m_{k}}\hookrightarrow E_{X}\right) =\inf \{\left\|
P\right\| :\text{ \ }P:E_{X}\rightarrow Y^{m_{k}}\text{; \ }P^{2}=P\}
\end{equation*}

tend to infinity.

Certainly, $E_{X}$ contains a subspace $l_{2}$ (since $E_{X}$ is an
approximate envelope and thanks to the Dvoretzki's theorem).

Clearly, this space has no complement in $E_{X}$: in a contrary case any
sequence of finite dimensional subspaces of $E_{X}$ that are uniformly
isomorphic to Euclidean ones should have uniformly bounded projection
constants.

Assume that $\left( E_{X}\right) ^{\ast }$ contains a subspace $Z$, which is
isomorphic to $l_{2}$. Then, according to [19], $Z$ is a complemented
subspace of $\left( E_{X}\right) ^{\ast }$ since it is of type 2. Hence, by
superreflexivity, $\left( E_{X}\right) ^{\ast \ast }=E_{X}$ contains a
complemented subspace $Z^{\ast }$, which is also isomorphic to the Hilbert
space. However, as was shown before, it is impossible.
\end{proof}

\begin{remark}
In the same way (i.e. using the presence in certain spaces of kind $E_{X}$
of an uncomplememted Hilbertian subspace) it may be obtained the following
extention of the theorem 12.
\end{remark}

\begin{theorem}
Let $\ast \ast \left( X^{f}\right) $ be a superreflexive class such that $%
S\left( \ast \ast \left( X^{f}\right) \right) =\{2\}$. Let $E_{X}$ be the
corresponding space of almost universal disposition. If $E_{X}$ is not of
type 2 then it fails to have the approximation property. 
\end{theorem}

\begin{problem}
Whether superreflexive spaces of kind $E_{X}$ either are isomorphic to the
Hilbert space or has no approximation property?
\end{problem}

\section{References}

\begin{enumerate}
\item  Lindenstrauss J., Rosenthal H.P. \textit{Authomorphisms in }$c_{o}$%
\textit{, }$l_{1}$\textit{\ and }$m$, Israel J. Math. \textbf{7} (1969)
227-239

\item  Enflo P. \textit{A counterexample to the approximation problem in
Banach spaces}, Acta Math. \textbf{130 }(1973) 309-317

\item  Bessaga C., Pe\l czy\'{n}ski A. \textit{Banach spaces non -
isomorphic to their Cartesian squares}, Bull Acad. Pol. Sci. Ser Math.
Astron. et Phys. \textbf{8} (1960) 77-80

\item  Semadeni Z. \textit{Banach spaces non-isomorphic to their Cartesian
squares}, Bull Acad. Pol. Sci. Ser. Math. Astron. et Phys.\textbf{\ 8}
(1960) 81-84

\item  Gowers W.T. \textit{Recent results in the theory of
infinite-dimensional Banach spaces}, Proc. of Int. Congress Math. Zurich, v. 
\textbf{2} (1994) 931-942

\item  Schwartz L. \textit{Geometry and probability in Banach spaces}, Bull.
AMS \textbf{4:2} (1981) 135-141

\item  Gurarii V.I.\textit{\ Spaces of universal disposition, isotropic
spaces and the Mazur problem on rotations in Banach spaces}, Sibirsk. Mat.
Journ. (in Russian) \textbf{7} (1966) 1002-1013

\item  Figiel T. \textit{An example of infinite dimensional reflexive Banach
space non-isomorphic to its Cartesian square}, Studia Math. \textbf{42:8}
(1972) 295 - 306

\item  Pe\l czy\'{n}ski A.\textit{\ Projections in certain Banach spaces},
Studia Math. \textbf{19} (1960) 209-228

\item  Stern J. \textit{Ultrapowers and local properties in Banach spaces},
Trans. AMS \textbf{240} (1978) 231-252

\item  Tokarev E.V. \textit{Injective Banach spaces in the finite
equivalence classes }(transl. from Russian), Ukrainian Mathematical Journal 
\textbf{39:6} (1987) 614-619

\item  Maurey B., Pisier G. \textit{S\'{e}ries de variables al\'{e}atoires
vectorielles ind\'{e}pendantes et propri\'{e}t\'{e}s g\'{e}om\'{e}triques
des espaces de Banach,} Studia Math. \textbf{58} (1976) 45-90

\item  Lindenstrauss J., Pe\l czy\'{n}ski A. \textit{Absolutely summing
operators in }$\mathcal{L}_{p}$\textit{-spaces and their applications},
Studia Math. \textbf{29 }(1968) 275-326

\item  Rosenthal H.P. \textit{Projections onto translation-invariant
subspaces of }$L_{p}(G)$, Memoirs AMS \textbf{63} (1968) 1-84

\item  Bennet G., Dor L.E., Gudman V., Johnson W.B., Newmann C.M. \textit{On
uncomplemented subspaces of }$L_{p}$, $1<p<2$, Israel J. Math. \textbf{26}
(1977) 178-187

\item  Tsirelson B.S. \textit{Not every Banach space contains }$l_{p}$%
\textit{\ or }$c_{0}$ (in Russian), Funct. Analysis and its Appl. \textbf{8:2%
} (1974) 57-60

\item  Kwapie\'{n} S. \textit{Isomorphic characterization of inner product
spaces by orthogonal series with vector valued coefficients}, Studia Math. 
\textbf{44} (1972) 583-595

\item  Pisier G. \textit{On the dimension of the }$l_{p}(n)$\textit{%
-subspaces of Banach spaces for }$1<p<2$, Trans. AMS \textbf{276} (1983)
201-211

\item  Maurey B. \textit{Quelques probl\`{e}mes de factorization d'op\'{e}%
rateurs lin\'{e}aires}, Actes du Congress Int. Math. Vancouver, 1974 v. 
\textbf{2} (1975) 75-79
\end{enumerate}

\end{document}